\documentclass[10pt,twoside]{amsart}

\newcommand{\bb}[1]{\mbox{\boldmath ${#1}$}}
\newcommand{\myeqref}[1]{(\ref{#1})}

\def\pn{\par\noindent}
\def\cen{\centerline}

\begin{document}

\leftline{\scriptsize \it J. Appl. Math. {$\&$} Computing Vol. {\bf\rm 16} (2004),
No. 1-2, pp. 509-513  }

\vspace{1.3cm}

\title{SOURCES OF HIGH LEVERAGE \\IN LINEAR REGRESSION MODEL}

\author{Myung Geun Kim}

\thanks{ {\scriptsize Received January 13, 2004.  Revised, March 9, 2004.  }
\newline\indent{\scriptsize $\copyright$ 2004 Korean Society for Computational
{$\&$} Applied Mathematics and Korean SIGCAM}}

\maketitle

\begin{abstract}
 Some reasons for high leverage are analytically investigated by decomposing leverage into
 meaningful components.  The results in this work can be used for remedial action
as a next step of data analysis.

\vskip 0.4cm


\end{abstract}

\pagestyle{myheadings}
\markboth{\centerline {\scriptsize M.G. Kim}}
         {\centerline {\scriptsize Sources of High Leverage }}

\bigskip
\bigskip
\bigskip

\cen{\bf 1.  Introduction}
\bigskip
\medskip

Inclusion of high leverage points in regression data can mislead our conclusion and cause some other
statistical problems.
Regression outliers may not be identified by looking at the least squares residuals
when the outliers are high leverage points  because high leverage points tend to 
have very small  residuals as the least sqares fit is pulled too much in the direction of
these outlying points ([1], [4], [6]). Gunst [2] and Mason and Gunst [5] argued that if a high leverage point
possesses extreme values on some regresors, then it can induce a collinearity among regressors.
Mason and Gunst[5] showed that collinearity can be increased without bound by increasing the
leverage of a point.
A high leverage point can hide or create collinearity ([1], [3], [7]).

Chatterjee and Hadi [1] gave some conditions for high leverage. However, no attempt has been made to uncover 
sources of high leverage. In this work we will make an analytic investigation of reasons for high leverage.
To this end two decompositions of leverage into meaningful components are derived. For these derivations
some preliminary results are also obtained. The results in this work can be used for remedial action
as a next step of data analysis.

\bigskip
\bigskip
\bigskip

\pagebreak

\cen{\bf 2.  Preliminaries}
\bigskip
\medskip

Consider a multiple linear regression model defined by
$$
{\bb y} = \beta_{0}{\bb 1}_{n} + {\bb X} {\bb \beta} +  {\bb \varepsilon},
$$
where ${\bb y}$ is an $n \times 1$ vector of observations on a response variable,
${\bb X}=({\bb x}_{1},...,{\bb x}_{n})^{T}$ is an $n \times p$ matrix of measurements on $p$ regressors,
${\bb 1}_{n}$ is the $n \times 1$ vector of all elements equal to one,
$\beta_{0}$ and ${\bb \beta}=(\beta_{1},...,\beta_{p})^{T}$ are unknown regression
coefficients, and ${\bb \varepsilon}$ is an $n \times 1$ vector of unobservable random errors.

Let ${\bb Z} = ({\bb 1}_{n}~~{\bb X})$. Then the hat matrix is
${\bb H} = {\bb Z}({\bb Z}^{T}{\bb Z})^{-1}{\bb Z}^{T}$. We have the following identity
([1])
$$
 {\bb H} = \frac{1}{n}{\bb 1}_{n}{\bb 1}_{n}^{T}
+ {\tilde {\bb X}}({\tilde {\bb X}}^{T}{\tilde {\bb X}})^{-1}{\tilde {\bb X}}^{T},
$$
where ${\tilde {\bb X}} = ({\bb I}_{n} - (1/n){\bb 1}_{n}{\bb 1}_{n}^{T}){\bb X}$ is
the centered  data matrix and ${\bb I}_{n}$ is the identity matrix of order $n$.
The leverage of the $r$-th observation is the $r$-th diagonal element $h_{rr}$ of ${\bb H}$
given by
\begin{equation}\label{eq:lev1}
h_{rr} = \frac{1}{n}(1 + D_{r}^{2})
\end{equation}
where
$D_{r}^{2}=n({\bb x}_{r}-\bar{\bb x})^{T}({\tilde {\bb X}}^{T}{\tilde {\bb X}})^{-1}
({\bb x}_{r}-\bar{\bb x})$ is the squared Mahalanobis distance from ${\bb x}_{r}$
to the data mean $\bar{\bb x} = (1/n)\sum_{j=1}^{n}{\bb x}_{j}$.
Eq. \myeqref{eq:lev1} implies that sources of high leverage can be investigated by
figuring out reasons for large $D_{r}^{2}$.

\bigskip
\bigskip
\bigskip

{\bf 2.1.  Regression of a regressor on the other regressors}
\bigskip
\medskip

The $i$-th column of ${\tilde {\bb X}}$ is written as ${\tilde {\bb x}}_{i}$.
Let ${\bb A}$ be the lower triangular matrix with positive diagonal elements
such that $(1/n){\tilde {\bb X}}^{T}{\tilde {\bb X}}={\bb A}{\bb A}^{T}$.
That is, ${\bb A}$ is the Cholesky root of the data covariance matrix.
We put ${\bb B}^{T} = {\bb A}^{-1}$ and denote the last column of ${\bb B}$
by ${\bb b} = (b_{1p},...,b_{pp})^{T}$.
We partition
\begin{equation}\label{eq:chol}
{\bb A} = \left[ \begin{array}{cc}
{\bb A}_{1} & {\bb 0} \\
{\bb a}_{1}^{T} & a_{pp}
\end{array} \right]
\end{equation}
such that ${\bb A}_{1}$ is the leading principal submatrix, having order $p-1$, of ${\bb A}$.
If we partition ${\tilde {\bb X}} = ({\tilde {\bb X}}_{1}~~{\tilde {\bb x}}_{p})$,
then the least squares estimator of the vector of regression coefficients
for the regression of ${\tilde {\bb x}}_{p}$ on the other regressors is easily computed as
\begin{eqnarray*}
 ({\tilde {\bb X}}_{1}^{T}{\tilde {\bb X}}_{1})^{-1}{\tilde {\bb X}}_{1}^{T}{\tilde {\bb x}}_{p}
&=& ({\bb A}_{1}{\bb A}_{1}^{T})^{-1}{\bb A}_{1}{\bb a}_{1} \cr
&=& -\frac{1}{b_{pp}}(b_{1p},...,b_{p-1,p})^{T}.
\end{eqnarray*}
The $r$-th residual from the regression of the $p$-th regressor on the other regressors
becomes
\begin{equation}\label{eq:resid1}
\frac{1}{b_{pp}}{\bb b}^{T}({\bb x}_{r}-\bar{\bb x}).
\end{equation}
The error sum of squares is computed as
$\mbox{SSE} = nb_{pp}^{-2}.$
Hence the $r$-th standardized residual is  ${\bb b}^{T}({\bb x}_{r}-\bar{\bb x})$
multiplied by a constant.

In general, similar results can be obtained for the regression of the $i$-th $(1 \leq i \leq p)$
regressor on the remaining regressors by appropriately permuting regressors.

\bigskip
\bigskip
\bigskip

{\bf 2.2.  Multiple correlation coefficient}
\bigskip
\medskip

Let $r_{i}^{2}$ be the squared multiple correlation coefficient of the $i$-th
regressor with the remaining regressors and
$s_{i}^{2} = (1/n)\sum_{r=1}^{n}(x_{ri}-\bar{x}_{i})^{2}$ be the data variance of the $i$-th
regressor,
where  $x_{ri}$ is the $i$-th element of ${\bb x}_{r}$ and
$\bar{x}_{i}$ is the $i$-th element of $\bar{\bb x}$.
Then the squared multiple correlation coefficient of the last regressor ${\tilde {\bb x}}_{p}$
with the other regressors is easily computed as
\begin{eqnarray} 
r_{p}^{2} &=& \frac{
{\tilde {\bb x}}_{p}^{T}{\tilde {\bb X}}_{1}({\tilde {\bb X}}_{1}^{T}{\tilde {\bb X}}_{1})^{-1}{\tilde {\bb X}}_{1}^{T}
{\tilde {\bb x}}_{p}}
{{\tilde {\bb x}}_{p}^{T}{\tilde {\bb x}}_{p}} \nonumber \cr
&=&  \frac{{\bb a}_{1}^{T}{\bb a}_{1}}
    {{\bb a}_{1}^{T}{\bb a}_{1}+a_{pp}^{2}}.
\end{eqnarray}
Further we have
\begin{equation}\label{eq:mcc2}
1-r_{p}^{2} = (s_{p}^{2}b_{pp}^{2})^{-1}
\end{equation}
since
$s_{p}^{2}=(1/n){\tilde {\bb x}}_{p}^{T}{\tilde {\bb x}}_{p}={\bb a}_{1}^{T}{\bb a}_{1}+a_{pp}^{2}$.

In general, we can compute $r_{i}^{2}$ $(1 \leq i \leq p)$ by appropriately permuting regressors.

\bigskip
\bigskip
\bigskip

\cen{\bf 3.  Some reasons for high leverage}
\bigskip
\medskip

In this section we provide two decompositions of leverage, using the results in Section 2, that
can explain some reasons for high leverage.
\bigskip
\medskip

{\bf 3.1.  Decomposition I}
\bigskip
\medskip

We denote a generic vector of $p$ regressor variables by  ${\bb x}$. Let ${\bb x}_{(i)}$ be
a vector obtained by interchanging the $i$-th and $p$-th regressors of
${\bb x}$ and ${\bb P}_{(i)}$ be the associated permutation matrix such that
${\bb x}_{(i)} = {\bb P}_{(i)}{\bb x}$.
Let ${\tilde {\bb X}}_{(i)} = {\tilde {\bb X}}{\bb P}_{(i)}$.
It is understood that ${\bb x}_{(p)} = {\bb x}$, ${\tilde {\bb X}}_{(p)} = {\tilde {\bb X}}$
and ${\bb P}_{(p)} = {\bb I}_{p}$.
Let ${\bb A}_{(i)}$ be the Cholesky root of
the data covariance matrix for ${\bb x}_{(i)}$, that is,
$(1/n){\tilde {\bb X}}_{(i)}^{T}{\tilde {\bb X}}_{(i)}={\bb A}_{(i)}{\bb A}_{(i)}^{T}$.
This  can be expressed as
$\{(1/n){\tilde {\bb X}}_{(i)}^{T}{\tilde {\bb X}}_{(i)}\}{\bb B}_{(i)} ={\bb A}_{(i)}$,
where ${\bb B}_{(i)}^{T} = {\bb A}_{(i)}^{-1}$. Comparison of the last columns of both sides
of the previous equation yields
$\{(1/n){\tilde {\bb X}}^{T}{\tilde {\bb X}}\}{\bb P}_{(i)}{\bb b}_{(i)}={\bb P}_{(i)}{\bb a}_{(i)}$ $(i=1,...,p)$,
where ${\bb a}_{(i)}$ and ${\bb b}_{(i)}$ are the last columns of
${\bb A}_{(i)}$ and ${\bb B}_{(i)}$, respectively. Collection of these $p$ equations into a matrix form
gives the inverse of the data covariance matrix as
\begin{equation}\label{eq:si}
\left(\frac{1}{n}{\tilde {\bb X}}^{T}{\tilde {\bb X}}\right)^{-1}
= [b_{(1)pp}{\bb P}_{(1)}{\bb b}_{(1)},\ldots,b_{(p)pp}{\bb P}_{(p)}{\bb b}_{(p)}],
\end{equation}
where $b_{(i)pp}$ is the last element of ${\bb b}_{(i)}$.
Inserting  \myeqref{eq:si} into $D_{r}^{2}$ with use of \myeqref{eq:mcc2} yields
 the following
decomposition of $D_{r}^{2}$
\begin{equation}\label{eq:dec1}
D_{r}^{2}=\sum_{i=1}^{p}(1-r_{i}^{2})^{-1/2}
[{\bb b}_{(i)}^{T}{\bb P}_{(i)}({\bb x}_{r}-\bar{\bb x})]
\left(\frac{x_{ri}-\bar{x}_{i}}{s_{i}}\right).
\end{equation}
Note that ${\bb A}_{(p)}={\bb A}$, ${\bb B}_{(p)}={\bb B}$, $b_{(p)pp} = b_{pp}$
and ${\bb b}_{(p)}={\bb b}$.

Together with \myeqref{eq:lev1}, eq. \myeqref{eq:dec1}  
reveals some sources of high leverage.
The $i$-th term of \myeqref{eq:dec1} shows that
the contribution of the $i$-th term to the
$r$-th leverage $h_{rr}$ depends on three components.
The first component $(1-r_{i}^{2})^{-1/2}$  will be large
whenever there is a high relationship between the $i$-th regressor and the set of
remaining regressors. In this case the contribution of $(1-r_{i}^{2})^{-1/2}$ to the
$r$-th leverage $h_{rr}$ will be  effective. The second component
${\bb b}_{(i)}^{T}{\bb P}_{(i)}({\bb x}_{r}-\bar{\bb x})$ is the standardized residual
 from the regression of the $i$-th regressor on the remaining regressors in the light of 
\myeqref{eq:resid1}, and it will be large in its absolute value
whenever the value of the $i$-th regressor $x_{ri}$ is far from the hyperplane formed by
the remaining regressors, that is, whenever $x_{ri}$ is an outlier for the regression
of the $i$-th regressor on the remaining regressors.
The third component $(x_{ri}-\bar{x}_{i})/s_{i}$ indicates a marginally standardized deviation
and it will be large in its absolute value whenever $x_{ri}$ is a marginally outlier.

\bigskip
\bigskip
\bigskip

{\bf 3.2.  Decomposition II}
\bigskip
\medskip

If we partition the data covariance matrix $(1/n){\tilde {\bb X}}^{T}{\tilde {\bb X}}  = {\bb A}{\bb A}^{T}$
according to the partition of ${\bb A}$ in \myeqref{eq:chol},
then the inverse of the resulting partitioned matrix is computed as
\begin{equation}\label{eq:decc0}
\left(\frac{1}{n}{\tilde {\bb X}}^{T}{\tilde {\bb X}}\right)^{-1} =
\left[ \begin{array}{cc}
({\bb A}_{1}{\bb A}_{1}^{T} )^{-1}+{\bb b}_{1}{\bb b}_{1}^{T} & b_{pp}{\bb b}_{1} \\
b_{pp}{\bb b}_{1}^{T}  &  b_{pp}^{2}
\end{array} \right],
\end{equation}
where ${\bb b}_{1} = (b_{1p},\ldots,b_{p-1,p})^{T}$.
Note that ${\bb A}_{1}{\bb A}_{1}^{T}$ is the data covariance matrix for the first
$p-1$ regressor variables.
Let ${\bb x}_{r(-p)} = (x_{r1},\ldots,x_{r,p-1})^{T}$
and $\bar{\bb x}_{(-p)} = (\bar{x}_{1},\ldots,\bar{x}_{p-1})^{T}$.
We write the squared Mahalanobis distance based on the first $p-1$ regressor variables as
$D_{r(-p)}^{2} = ({\bb x}_{r(-p)}-\bar{\bb x}_{(-p)})^{T}({\bb A}_{1}{\bb A}_{1}^{T} )^{-1}
({\bb x}_{r(-p)}-\bar{\bb x}_{(-p)})$. Then a little computation
with use of \myeqref{eq:decc0}
 gives the following decomposition of
$D_{r}^{2}$
\begin{equation}\label{eq:decc1}
D_{r}^{2} = D_{r(-p)}^{2} + [{\bb b}^{T}({\bb x}_{r}-\bar{\bb x})]^{2}.
\end{equation}
In view of this decomposition,
the removal of the $p$-th regressor from the regression model  decreases the value of $h_{rr}$
  by the second term of \myeqref{eq:decc1} divided by $n$
  when the remaining regressors are still kept in the model.

\end{document}